\documentclass[12pt]{article}
\usepackage[a4paper, left=2cm, right=2cm, bottom=3cm, top=3cm]{geometry}
\usepackage{amsmath}
\usepackage{srcltx}
\usepackage{amsfonts}
\usepackage{amssymb}
\usepackage{lipsum}
\usepackage{psfrag}
\usepackage{verbatim}
\usepackage{graphicx} 
\usepackage{color}
\usepackage{faktor}
\providecommand{\U}[1]{\protect\rule{.1in}{.1in}}

\newtheorem{theorem}{Theorem}

\newtheorem{claim}[theorem]{Claim}

\newtheorem{definition}[theorem]{Definition}
\newtheorem{example}[theorem]{Example}

\newtheorem{lemma}[theorem]{Lemma}

\newtheorem{proposition}[theorem]{Proposition}
\newtheorem{remark}[theorem]{Remark}

\newcommand*{\QEDA}{\hfill\ensuremath{\blacksquare}}

\newcommand{\adj}{\mathop{\rm ad}\nolimits}

\newcommand{\R}{\ensuremath{\mathbb{R}}}

\newcommand{\N}{\ensuremath{\mathbb{N}}}

\newcommand{\Z}{\ensuremath{\mathbb{Z}}}

\begin{document}
\title{Controllability of discrete-time linear systems on Lie groups with finite semisimple center}
\author{Thiago M. Cavalheiro\\Departamento de Matem\'{a}tica, Universidade Estadual de Maring\'{a}\\Maring\'{a}, Brazil.
\and João A. N. Cossich\\Departamento de Matem\'{a}tica, Universidad de Tarapac\'{a} - Arica, Chile
\and Alexandre J. Santana\\Departamento de Matem\'{a}tica, Universidade Estadual de Maring\'{a}\\Maring\'{a}, Brazil }
\maketitle

\begin{abstract} In this paper we stated a condition for the controllability of discrete-time linear systems for the case when the Lie group has finite semisimple center.
\end{abstract}



\section{Introduction} 

The aim of this paper is to study the controllability of a discrete-time linear control system in the form 
\begin{equation*}
    \Sigma: x_{k+1} = f_{u_k}(x_k),  
\end{equation*}
with restricted control $u_k \in U \subset \mathbb{R}^m$ a compact convex neighborhood of $0$, under the hypothesis that the state space $G$ is a connected Lie group with finite semisimple center. Despite having few tools to deal with such systems, the study of the eigenvalues of the function $(df_0)_e$ have been showing a accurate way to prove some consistent properties of the reachable sets. For instance, in the case when $G$ is a solvable connected Lie group and the automorphism $f_0: G \longrightarrow G$ of the system is inner, Cavalheiro, Santana and Cossich \cite{TAJ} proved that a sufficient condition for controllability is that $(df_0)_e$ has only eigenvalues with modulus 1 and the reachable set of the identity is open. When $G$ is nilpotent, this condition is also necessary. This paper is particularly interesting given that the discretization of the continuous case (see Ayala and da Silva \cite{AyalaeAdriano2}) is a particular example of the present case, as we will discuss in the last part of this paper. 

Regarding the continuous case, in Jouan \cite{jouan} it is stated that the flow $\varphi_t$ of the linear vector field is associated with a derivation $\mathcal{D}$ on the Lie algebra $\mathfrak{g}$ in the following way $e^{t\mathcal{D}} = d\varphi_t$. In our type of system, the automorphism $f_0 \in \hbox{Aut}(G)$ have no reason at all to be as described. In this case, we are considering the case when $f_0$ is a inner automorphism in the form $df_0 = e^{\mathcal{D}_1} \circ ... \circ e^{\mathcal{D}_k}$, for some $k \geq 1$. Under some hypothesis of $df_0$, denoting by $\mathcal{R}$ the reachable set of $e \in G$, we prove the following statement, which is the main result of this paper: \textit{Let $G$ be a connected Lie subgroup with finite semisimple center. If $e \in \hbox{int}\mathcal{R}$, $\mathcal{D}_i \circ \mathcal{D}_j =\mathcal{D}_j \circ \mathcal{D}_i $ for every $i,j=1,...,k$ and every eigenvalue of $df_0$ has norm 1, the system $(\Sigma)$ is controllable.} This paper ends up with some examples in the semisimple Lie group $\hbox{SL}_2(\R)$. 

This paper is structured as follows: the section $2$ is dedicated for some basic concepts of control systems and useful results used along this paper. In section $3$ we proved the main results using an auxiliary system and at the last part of this section, we explored the case when the function $f_0$ is a inner automorphism of $G$. The section $4$ is dedicated for the construction of the class of discrete-time linear systems on $\hbox{SL}_2(\R)$ and some particular examples.

\section{Preliminaries}

\subsection{General properties} This section will be dedicated to defined the general properties of control systems and show the terminology that will be used along this paper. Initially, our phase space  $\mathcal{C}^{\infty}-$ is a Riemannian manifold $n-$dimensional $M$ endowed with a canonical metric $d$ and taking a $\mathcal{C}^{\infty}$ map $f: U \times M \longrightarrow M$, defined over $U$ a non-empty compact convex neighborhood of $0 \in \mathbb{R}^m$ such that $U \subset \overline{\hbox{int} U}$, using the notation $f_u: M \longrightarrow M$ for $f_u(x) = f(u,x)$, the system we will study is in the form
\begin{equation}\label{sistgeneral}
    \Sigma: x_{k+1} = f_{u_k}(x_k), u_k \in U,
\end{equation}
with $k \in \mathbb{N}_0 = \N \cup \{0\}.$ Also, we will consider that for a open set $\hat{U}$ containing $U$, the function $f: \hat{U} \times M \longrightarrow M$ is a $\mathcal{C}^{\infty}$ map. 

Given any $x \in M$ as initial condition, the solution of $(\Sigma)$ will be denoted by $\varphi(k,x,u)$, where $u \in \mathcal{U} = \prod_{i \in \Z} U$ such that $u = (u_i)_{i \in \mathbb{Z}}$. For any $x \in M$, if for each $u \in U$ the function $f_u: M \longrightarrow M$ is a homeomorphism, the solution $\varphi$ is also defined for $k \in \Z$. Assuming such property, the solution has the form
\begin{equation*}
    \varphi(k,x_0,u) = 
    \left\{
    \begin{array}{cc}   
        f_{u_{k-1}} \circ ... \circ f_{u_0}(x_0),& k > 0\\
        x_0,& k=0\\
        f_{u_{k}}^{-1} \circ ... \circ f_{u_{-1}}^{-1}(x_0),& k < 0
    \end{array}
    \right.
\end{equation*}
and considering the function $\Theta: \Z \times \mathcal{U} \longrightarrow \mathcal{U}$, defined by $\Theta_k((u_i)_{i \in \Z}) = (u_{i+k})_{i \in \Z}$, the solution $\varphi$ satisfies the cocycles property which means that
\begin{equation*}
    \varphi(k+t, x,u) = \varphi(k,\varphi(t,x,u), \Theta_t(u))= \varphi(t,\varphi(k,x,u), \Theta_k(u)), \forall k,t \in \Z.  
\end{equation*}

The solution $\varphi$ also satisfy the following property: if $ts > 0$ in $\Z$, given $u,v \in \mathcal{U}$, there is a $w \in \mathcal{U}$ such that 
\begin{equation*}
    \varphi(t,\varphi(s,g,u),v) = \varphi(t + s, g, w), \forall g \in M.
\end{equation*}

All spaces will be considered endowed with the canonical topology.  Therefore, the space shift space $\mathcal{U}$ is compact.

\begin{remark} The function $\Theta$ is continuous and also satisfies the properties $\Theta_{t + s}(u) = \Theta_t(\Theta_s(u))$ and $\Theta_0(u) = u$. Then, $\Theta$ defines on $\mathcal{U}$ a continuous dynamical system. 
\end{remark}

\begin{definition}For $x \in M$, the set of points reachable and controllable from $x$ up to time $k > 0$ in $\mathbb{N}$ are
\begin{eqnarray*}
    \mathcal{R}_k(x) = \{y \in M: \hbox{ there is }u \in \mathcal{U} \hbox{ with }\varphi(k,x,u)=  y\}\\
    \mathcal{C}_k(x) = \{y \in M: \hbox{ there is }u \in \mathcal{U} \hbox{ with }\varphi(k,y,u)= x\}\\
\end{eqnarray*}
The sets $\mathcal{R}(x) = \bigcup_{k \in \mathbb{N}} \mathcal{R}_k(x)$ and $\mathcal{C}(x) = \bigcup_{k \in \mathbb{N}} \mathcal{C}_k(x)$ denote the reachable set and the controllable set from $x$ respectively.  
\end{definition}

\begin{definition}\label{regular}For each $k \in \N$, consider the function $G_k(g,u) = \varphi(k,g,u)$.  A pair $(g,u) \in M \times \hbox{int}U^k$ is called regular if $\hbox{rank}\left[\frac{\partial}{\partial u} G_k(g,u)\right] = \hbox{dim}M.$ We denote by 
\begin{equation*}
    \hat{\mathcal{R}}_k(g) = \{\varphi(k,g,u): (g,u) \in M\times \hbox{int}U^k \hbox{ is regular}\}.
\end{equation*}
the regular trajectory of $g \in G$ up to time $k \in \N$ and the regular trajectory  of $g \in G$ by $\hat{\mathcal{R}}(g) = \bigcup_{k \in \N} \hat{\mathcal{R}}_k(g)$. 
\end{definition}

In particular, the set $\hat{\mathcal{R}}(g)$ is open, for every $g \in M$. The system (\ref{sistgeneral}) is forward accessible (resp. backward accessible) if $\hbox{int}\mathcal{R}(x) \neq \emptyset$ (resp. $\hbox{int}\mathcal{C}(x) \neq \emptyset$), for all $x \in M$. The system $(\ref{sistgeneral})$ is accessible if both conditions are satisfied. Let us consider $M = G$ as a connected $n-$dimensional Lie group. 

\begin{definition}\label{linearsystem} A discrete-time control system
\begin{equation*}
    \Sigma: x_{k+1} = f_{u_k}(x_k), u_k \in U,
\end{equation*}
defined over $G$ with $U \subset \R^m$ a compact neighborhood of $0$ is said to be linear if 
\begin{itemize}
    \item[1-] $f_0: G \longrightarrow G$ is an automorphism; 
    \item[2-] The function $f$ satisfies 
        \begin{equation}\label{transl}
            f_u(g) = f_u(e) \cdot f_0(g). 
        \end{equation}
    where $"\cdot"$ denotes de product on $G$.
\end{itemize}
\end{definition}

For linear systems, the function $f$ can be defined using the translations of $G$. Given $u \in U$, as $f_u(e) \in G$, the expression (\ref{transl}) allow us to write $f_u(g)$ as 
\begin{equation}\label{transl1}
    f_u(g) = f_u(e)f_0(g) = L_{f_u(e)}(f_0(g)), 
\end{equation}
where $L_{f_u(e)}$ is the left translation by the element $f_u(e)$. Considering the expression above, the inverse of $f_u$ is given by 
\begin{equation}
    (f_u)^{-1}(g) = f_0^{-1}\circ L_{(f_u(e))^{-1}}(g)= f_0^{-1}((f_u(e))^{-1} \cdot g). 
\end{equation}

Then, we can conclude that $f_u$ is a diffeomorphism of $G$, for any $u \in U$. The solutions can also be defined in terms of translations of the solution at the neutral element, as in the next proposition. The proof can be found in (Colonius, Cossich and Santana \cite{CCS1}). 

\begin{proposition}\label{prop52} Consider a discrete-time linear control system $x_{k+1} = f(u_k,x_k)$, $u_k \in U$ defined on a Lie group $G$. Then it follows for all $g \in G$ and $u = (u_i)_{i \in \Z} \in \mathcal{U}$ that
\begin{equation*}
    \varphi(k,g,u) = \varphi(k,e,u)f_0^k(g). 
\end{equation*}
\end{proposition}

By using the reversed-time system 
\begin{equation*}
    \Sigma': g_{k+1} = \hat{f}_{u_k}(g_k), k \in \N_0,
\end{equation*}
given by the function $\hat{f}_u(g) = f_u^{-1}(e)f_0^{-1}(g)$, denoting by $\mathcal{R}^*_k$ and $\mathcal{C}^*_k$ the reachable and controllable sets up to time $k$ of $e$ of the system $(\Sigma')$, its is proved in \cite{TAJ} the following result. 

\begin{lemma}\label{setsinvsys}It holds that $\mathcal{R}_k^*=\mathcal{C}_k$ and $\mathcal{R}_k=\mathcal{C}_k^*$, for all $k\in\mathbb{N}$.
\end{lemma}

\section{Conditions for controllability}

Let us consider a connected Lie group $G$ with Lie algebra $\mathfrak{g}$ and the discrete-time linear system 
\begin{equation*}
    \Sigma: g_{k+1} = f_{u_k}(g_k), n \in \N_0,
\end{equation*}
with $u_k \in U$ compact convex neighborhood of $0 \in \R^m$. As said before, the function $f_u: G \longrightarrow G$ is a diffeomorphism for any $u \in U$ and $f_0$ is a automorphism of $G$. Then the system $\Sigma$ is defined for any $k \in \Z$. Denote by $\mathfrak{g}$ the Lie algebra of $G$, endowed with the Lie bracket 
\begin{equation}\label{liebracket}
    [X,Y] = \frac{\partial^2}{\partial t \partial s} \left(X_{-t} \circ Y_s \circ X_t\right)\bigg|_{t=s=0}. 
\end{equation}
where $X_t$ and $Y_t$ are the respectives exponentials of $X$ and $Y$ at the time $t \in \R$. It is very-known (see \cite{sanmartin1}) that $[X,Y] = 0$ if, and only if, $\exp{X}\exp{Y} = \exp{Y}\exp{X}$. 

Considering the reachable sets $\mathcal{R}_k(e)$, $\mathcal{R}_{\leq k}(e) = \{\varphi(t,e,u): t \in [0,k] \cap \N, u \in \mathcal{U}\}$ and $\mathcal{R}(e) = \bigcup_{k \in \N}\mathcal{R}_k(e)$, it is easy to see that $e \in \mathcal{R}_k(e)$ for any $k \in \N.$ Besides, using the notation $\mathcal{R}(e) = \mathcal{R}$, $\mathcal{R}_k(e) = \mathcal{R}_k$ and $\mathcal{R}_{\leq k} = \mathcal{R}_{\leq k}(e)$, we get the following properties of $\Sigma$. The next two results will be often used along this paper and its proof can be found in \cite{TAJ}.

\begin{proposition}\label{reachablesetprop}The reachable set $\mathcal{R}$ satisfy the following properties: 
\begin{itemize}
    \item[1-] Given $\tau \geq 1$ in $\N$, then $\mathcal{R}_{\tau} = \mathcal{R}_{\leq \tau}$. 
    \item[2-] Given $0 < \tau_1 \leq \tau_2$ in $\N$, then $\mathcal{R}_{\tau_1} \subset \mathcal{R}_{\tau_2}$. 
    \item[3-] If $g \in G$, then $\mathcal{R}_{\tau}(g) = \mathcal{R}_{\tau} f_0^{\tau}(g)$. 
    \item[4-] If $\tau_1, \tau_2 \in \N$, then $\mathcal{R}_{\tau_1 + \tau_2} = \mathcal{R}_{\tau_1} f_0^{\tau_1}(\mathcal{R}_{\tau_2}) = \mathcal{R}_{\tau_2} f_0^{\tau_2}(\mathcal{R}_{\tau_1})$. 
    \item[5-] For any $u \in \mathcal{U}$, $g \in G$ and $k \in \N$, then 
    \begin{equation*}
        \varphi(k,\mathcal{R}(g),u) \subset \mathcal{R}(g).
    \end{equation*}
    \item[6-] $e \in \hbox{int}\mathcal{R}$ if and only if $\mathcal{R}$ is open. 
\end{itemize}
\end{proposition}

\begin{lemma}\label{AginA} Let be $g \in \mathcal{R}$ and assume that $f_0^t(g) \in \mathcal{R}$ for any $t \in \Z$. Then $\mathcal{R}\cdot g \subset \mathcal{R}$. 
\end{lemma}

Considering the system $(\Sigma),$ regarding the set $\mathcal{R}$ and also denoting $\mathcal{C}_k(e) = \mathcal{C}_k$, we have the following connection between the sets $\mathcal{R}_k$ and $\mathcal{C}_k$. 

\begin{proposition}\label{openess}Consider the linear system $(\Sigma)$. Then $\hbox{int}\mathcal{R}_k \neq \emptyset$ if, and only if, $\hbox{int}\mathcal{C}_k \neq \emptyset.$  
\end{proposition}

\noindent\textit{Proof:} Let us suppose $\hbox{int}\mathcal{R}_k \neq \emptyset$ and consider the automorphism $f_0$ of $(\Sigma)$. Take $g \in \hbox{int}\mathcal{R}_k$. Hence, there is a $u \in \mathcal{U}$ such that 
\begin{equation*}
    g = \varphi(k,e,u),
\end{equation*}
that is $\varphi(k,e,u)g^{-1} = e$. Using the properties of $\varphi$ we get 
\begin{equation*}
    \varphi(k,e,u)g^{-1} = \varphi(k,e,u)f_0^k(f_0^{-k}(g^{-1})) = \varphi(k,f_0^{-k}(g^{-1}),u) = e,
\end{equation*}
that is $f_0^{-k}(g^{-1}) \in \mathcal{C}_k$. Now, consider $V$ a neighborhood of $g$ such that $g \in V \subset \mathcal{R}_k.$ By the arguments above, $f^{-k}_0(V^{-1}) \subset \mathcal{C}_k$. The function $f_0^{-k}$ is an automorphism of $G$ and $V^{-1}$ is a neighborhood of $g^{-1}$. Then $f_0^{-k}(g^{-1}) \in \hbox{int}\mathcal{C}_k$.

Now, if $\hbox{int}\mathcal{C}_k\neq \emptyset$, let us take $g \in \hbox{int}\mathcal{C}_k$. Then, there are $k \in \N$ and $u \in \mathcal{U}$ such that $\varphi(k,g,u) = e$. Hence
\begin{equation*}
    \varphi(k,g,u) = \varphi(k,e,u)f_0^k(g) = e,
\end{equation*}
and consequently $\varphi(k,e,u) = f_0^{k}(g^{-1})$. If $g \in V \subset \mathcal{C}_k$, for some neighborhood $V$ of $g$, by the previous argument, $f_0^k(V^{-1}) \subset \mathcal{R}_k$. Therefore $\hbox{int}\mathcal{R}_k \neq \emptyset.$ \QEDA

\begin{remark}For the case when $\hbox{int}\mathcal{R} \neq \emptyset$, in \cite{CCS2} it is cited that there is a $k_0 \geq 1$ such that $\hbox{int}\mathcal{R}_k \neq \emptyset$ for every $k \geq k_0$. In this case, the proposition above ensures that $\hbox{int}\mathcal{R} \neq \emptyset$ if, and only if, $\hbox{int}\mathcal{C} \neq \emptyset$. 
\end{remark}

Given an automorphism $f \in \hbox{Aut}(G)$, for the case when $\mathfrak{g}$ is defined over a closed field, using the notation $df = d(f)_e: \mathfrak{g} \longrightarrow \mathfrak{g}$, following the concepts of da Silva, Ayala and Roman-Flores \cite{AyalaeRomaneAdriano}, we can consider the generalized eigenspaces of $df$ by
\begin{equation*}
    \mathfrak{g}_{\alpha} = \{X \in \mathfrak{g}: (df - \alpha)^n X = 0, \hbox{ for some }n \in \N\}, 
\end{equation*}
associated with the eigenvalue $\alpha$, we get the sets 
\begin{equation}\label{generalizeigenspaces}
    \mathfrak{g}^+ = \bigoplus_{|\alpha|>1} \mathfrak{g}_{\alpha}, \hbox{ } \mathfrak{g}^- =\bigoplus_{|\alpha|<1} \mathfrak{g}_{\alpha}, \hbox{ } \mathfrak{g}^0 = \bigoplus_{|\alpha|=1} \mathfrak{g}_{\alpha}. 
\end{equation}
and the primary decomposition of $\mathfrak{g}$ 
\begin{equation}\label{decomp}
    \mathfrak{g} = \mathfrak{g}^+ \oplus \mathfrak{g}^0 \oplus \mathfrak{g}^-. 
\end{equation}

Such decomposition is also true for the real case (see \cite{AyalaeRomaneAdriano} for further explanation). Based on the sets above, let us consider the following Lie subgroups of $G$: $G^0=\langle \exp{\mathfrak{g}^0} \rangle$, $G^+ = \exp{\mathfrak{g}^+}$ and $G^- = \exp{\mathfrak{g}^-}$. In $G^+$ and $G^-$, since $\mathfrak{g}^+$ and $\mathfrak{g}^-$ are nilpotent subalgebras of $\mathfrak{g}$, the exponential function is surjective. Given a homomorphism $\phi$ of $G$, we say that a Lie subgroup $H$ of $G$ is $\phi-$invariant if $\phi(H) \subset H$. On the system $\Sigma$, as the function $f_0$ is a automorphism, if the Lie subgroup $H$ is $f_0-$invariant, is easy to see that $f_0(H) = H$. In this case, the invariance also satisfies $f_0^{-k}(H) \subset H$, for any $k \in \N$. From now on, we will suppose that $\mathcal{R}$ is a open set of $G$. The next lemma will often be used and can be found in (Wustner \cite[Lemma 3.1]{sontag1}).

\begin{lemma}\label{lema212}Let $G$ be a Lie group with Lie algebra $\mathfrak{g}$ and $N$ a normal Lie subgroup of $G$ with Lie algebra $\mathfrak{n}$. Then for every $X \in \mathfrak{g}$, we have that 
\begin{equation*}
    \exp{(X + \mathfrak{n})}\subset \exp{(X)}N. 
\end{equation*}
\end{lemma}

From now on in this section, we will prove some results associated with the controllability of the system $\Sigma$. The next results will be helpful to our purposes. Denoting by $\mathcal{R}$ and $\mathcal{C}$ the respectives reachable of and controllable sets of the neutral element of $\Sigma$, taking the primary decomposition (\ref{decomp}) of $f_0$, we have the following result. 

\begin{lemma}\label{lemmaimport} 
Let $N \subset G^0$ be a $f_0-$invariant connected solvable Lie subgroup of $G^0$. Then $G^0 \subset \mathcal{R}$.   
\end{lemma}

The case we are considering is when $G$ satisfy the following definition.  

\begin{definition}\label{finitecenter1} Let $G$ be a connected Lie group. We say that $G$ has finite semisimple center if every semisimple Lie subgroup of $G$ have finite center.     
\end{definition}

Consider the linear system
\begin{equation}\label{linsyst}
    \Sigma: g_{k+1} = f_{u_k}(g_k), k \in \N_0,
\end{equation}
defined by the function $f: U \times G \longrightarrow G$, with $u_k \in U$ a compact convex neighborhood of $0$ of $\R^m$. In this particular case, are considering that $df_0$ satisfies 
\begin{equation*}
    df_0(X) = e^{\mathcal{D}_1} \circ ... \circ  e^{\mathcal{D}_n}(X), X \in \mathfrak{g}, 
\end{equation*}
where $\mathcal{D}_j$ are derivations of $\mathfrak{g}$. In this paper, our primary objective is to demonstrate the following theorem.

\begin{theorem}Let $G$ be a Lie group with finite semisimple center. If every eigenvalue of $\mathcal{D}_i$ has zero real part, $\mathcal{D}_i \circ \mathcal{D}_j = \mathcal{D}_j \circ \mathcal{D}_i$ for every $j,i=1,...,n$ and $\mathcal{R}$ is open, then the system $(\Sigma)$ is controllable.   
\end{theorem}

Let us start with the case $df_0 = e^{\mathcal{D}}$, for some $\mathcal{D}: \mathfrak{g} \longrightarrow \mathfrak  {g}$ derivation of $\mathfrak{g}$. Considering the solvable radical $\mathfrak{r}(\mathfrak{g}^0)$ of $\mathfrak{g}^0$, we have that $\mathfrak{s} = \mathfrak{g}^0/\mathfrak{r}(\mathfrak{g}^0)$ is a semisimple Lie algebra. As the solvable radical is invariant by automorphisms, we have that $e^{\mathcal{D}} \mathfrak{r}(\mathfrak{g}^0) =  \mathfrak{r}(\mathfrak{g}^0).$ Hence $\mathfrak{r}(\mathfrak{g}^0)$ is $\mathcal{D}-$invariant. Considering the canonical projection $\pi^*: \mathfrak{g}^0 \longrightarrow  \mathfrak{s}$, this also implies that the derivation $\mathcal{D}^*: \mathfrak{s} \longrightarrow \mathfrak{s}$ defined by $\mathcal{D}^*(\Bar{X})= \pi^*(\mathcal{D}(X))$ is well-defined on $\mathfrak{s}$ and satisfies the commuting property $\pi^* \circ \mathcal{D} = \mathcal{D}^* \circ \pi^*$. Since any derivation in a semisimple Lie algebra is inner \cite[Proposition 3.14]{sanmartin2}, there is a $Z \in \mathfrak{g}^0$ such that $\mathcal{D}^* = \hbox{ad}_s(\pi(Z))$.  

Considering $Z \in \mathfrak{g}^0$ as before, it turns out that 
\begin{equation*}
    df_0^k X = e^{\hbox{ad}(Z)}X + Z, X \in \mathfrak{g}^0.   
\end{equation*}
for some $Z = Z_{k,X} \in \mathfrak{r}(\mathfrak{g}^0).$ As $\mathfrak{r}(\mathfrak{g}^0)$ is a ideal in $\mathfrak{g}^0$, by Proposition (\ref{lema212}) we get that 
\begin{equation*}
    \exp( df_0^k X) = \exp(e^{\hbox{ad}(Z)}X + Z) = \exp(e^{\hbox{ad}(Z)}X  )\cdot g, 
\end{equation*}
for some $g = g_{X,k} \in R$, that is, 
\begin{equation*}
    f_0^k(\exp(X)) = \exp(kZ) \exp(X) \exp(kZ)g, g = g_{X,k} \in R.  
\end{equation*}

By connectedness of $G^0$, we get $f_0(h) = \exp(kZ)h\exp(kZ)g, g = g_{h,k} \in R$. Now, considering the set $\mathcal{R}_k(\exp(kY))$, by the previous analysis we get 
\begin{equation*}
    \mathcal{R}_k(\exp(kY)) = \mathcal{R}_k f_0^k(\exp(kY)) = \mathcal{R}_k \exp(kY) g, g \in R.  
\end{equation*}

\subsection{Semigroups and Controllability}

Considering the canonical projection $\pi: G^0 \longrightarrow S$, we get based on the preceding logic that 
\begin{equation*}
    \pi(G^0 \cap \mathcal{R}_k(\exp(kY))) = \pi((G^0 \cap \mathcal{R}_k) \exp(kY)).
\end{equation*}

Let us take $S = G^0/R$, the semisimple Lie group with Lie algebra $\mathfrak{s}$. Take the sets $\mathcal{S}_k = \pi(G^0 \cap \mathcal{R}_k(e^{kZ}))$ and $\mathcal{S} = \bigcup_{k \in \N} \mathcal{S}_k$. In particular we have $\mathcal{S}_k = \pi(G^0 \cap \mathcal{R}_k(e^{kZ})) = \pi((G^0 \cap \mathcal{R}_k)e^{kZ})$. Then we have the following proposition. 

\begin{proposition}If $\mathcal{R}$ is open, the set $\mathcal{S}$ is a semigroup of $S$ with non-empty interior. 
    \end{proposition}

\noindent\textit{Proof: } In fact, taking $x_1, x_2 \in \mathcal{S}$, we have that there are $u,v \in \mathcal{U}$ and $k,s \in \N$ such that 
\begin{equation*}
    x_1 = \pi(\varphi(k,e,u)e^{kZ}), x_2 = \pi(\varphi(s,e,v)e^{sZ})
\end{equation*}
and $x_1x_2 = \pi(\varphi(k,e,u)e^{kZ}\varphi(s,e,v)e^{kZ})$. Considering $y=\varphi(s+k,e^{(s+k)Z},w))$, where $w$ is the concatenation between $v$ and $u$, using the fact of $\Sigma^C$ is linear we get 
\begin{eqnarray*}
    y &=& \varphi(k+s,e^{(k+s)Z},w) = \varphi(k,\varphi(s,e^{(k+s)Z},v),u) = \varphi(k,e,u)(h_0^C)^k(\varphi(s,e^{(k+s)Z},v)) \\
    &=&\varphi(k,e,u)e^{kZ}(\varphi(s,e^{(k+s)Z},v))e^{-kZ}  = \varphi(k,e,u)e^{kZ}\varphi(s,e,v)(h_0^C)^s(e^{(k+s)Z}) e^{-kZ}\\
    &=& \varphi(k,e,u)e^{kZ}\varphi(s,e,v)(e^{sZ} e^{(s+k)Z}e^{-sZ}) e^{-kZ} = \varphi(k,e,u)e^{kZ}\varphi(s,e,v)e^{sZ}. 
\end{eqnarray*}
that is 
\begin{equation*}
    \pi(y) = \pi(\varphi(k,e,u)e^{kZ}\varphi(s,e,v)(e^{sZ})) = x_1x_2 \in \mathcal{S}_{s+k}. 
\end{equation*}

Hence $x_1x_2 \in \pi((G^0 \cap \mathcal{R}_{k_1+k_2})e^{(k_1+k_2)Z}) = \mathcal{S}_{k_1+k_2} \subset \mathcal{S}$, that is, $\mathcal{S}$ is a semigroup of $S$. Now, we have by hypothesis that $e \in \hbox{int}\mathcal{R}.$ Then there is a $k_0 \in \N$ such that $e \in \hbox{int}\mathcal{R}_{k_0}$. Take $X = k_0 \pi(Z)$. We have that 
\begin{equation*}
    \exp_s{(X)}= \pi(\exp(k_0 Z)) \in \pi((\hbox{int}\mathcal{R}_{k_0} \cap G^0)e^{k_0Z}) \subset \hbox{int} \mathcal{S}. 
\end{equation*}

Then $\hbox{int}\mathcal{S} \neq \emptyset.$ \QEDA

Now, let us prove the main result associated with the system $\Sigma$. 

\begin{proposition}\label{h_C}Let $G$ be a connected Lie group with finite semisimple center. If $e \in \hbox{int}\mathcal{R}$, then $G^0 \subset \mathcal{R}$. 
\end{proposition}

\noindent\textit{Proof:} Without loss of generality, let us consider $\mathfrak{g}$ as the right-invariant vector fields of $G$. For $S = G^0/R$, where $R= \langle \exp{\mathfrak{r}(\mathfrak{g}^0)} \rangle$ we have two possibilities for $S$: 

If $S$ is compact, using the fact that any semigroup with nonempty interior of a compact group contains the identity component of the group, and considering that $S$ is connected, we have $\mathcal{S} = S$. 

Now, if $S$ is noncompact, consider the result \cite[Lemma 4.1]{sanmartin3}: 

\begin{lemma}\label{lemma} Let $G$ be a connected Lie group with finite semisimple center and $\mathcal{S}$ a semigroup of $G$. If there exists $X \in \mathfrak{g}$ such that $\hbox{ad}(X)$ is nilpotent and $e^X \in \hbox{int}\mathcal{S}$, then $G = \mathcal{S}$.   
\end{lemma}

Following with the proof, for the $X$ in the previous result, let us consider the Jordan decomposition of $X = X_n + X_s$, with $\hbox{ad}X_n$ nilpotent and $\hbox{ad}X_s$ semisimple\footnote{Considering the complexification $\mathfrak{g}_{\mathbb{C}}$, the operator $\hbox{ad}X_n$ is diagonalizable.}. As $X \in \mathfrak{g}^0/\mathfrak{r}(\mathfrak{g}^0)$, then $e^{\hbox{ad}X} = e^{\hbox{ad}X_s}$ has the same spectrum of $d\Bar{f}_0$, with all eigenvalues with norm $1$. As $S$ has finite center, the set 
\begin{equation*}
    K = \hbox{cl}\{\exp_s(\tau X_s): \tau \in \R\}, 
\end{equation*}
is compact. As a matter of fact, the element $e^{X_s}$ is a compact element of $S$. To prove that, as $S$ is connected, we get $\ker\hbox{Ad} = Z(S)$. As $S$ is semisimple with finite center, $Z(S)$ is finite. As $e^{\hbox{ad}(X_s)}$ has all eigenvalues with norm $1$, there is a $A \subset \hbox{Aut}(\mathfrak{g})$ a compact subgroup containing $e^{\hbox{ad}(X_s)}$ \cite[Corollary 1.2]{fresnel}. Consider the set $B = \hbox{Ad}^{-1}(A)$. The set $B$ is compact, given that $B/Z(S) \cong A$ and $Z(S)$ is normal \cite[Page 61, item F]{pontr}. Also, $B$ contains $e^{X_s}.$ Then $e^{X_s}$ is a compact element. Now, consider the $1$-parameter subgroup $L = \{\exp_s(\tau X_s): \tau \in \R\}$. Then $L$ is closed or has compact closure \cite[Exercise c.2]{helgason}. If $L$ is closed and non-compact, then $L \cong \mathbb{R}$, which does not contains any non-trivial elements, contradicting $e^{X_s} \in L$. Then $K = \hbox{cl}L$ is compact. 

Take the subset 
\begin{equation*}
    K_s = \{x \in K: x \exp(\tau X_n) \in \hbox{int}\mathcal{S}, \hbox{for some }\tau \in \R\}. 
\end{equation*}

Then $\hbox{int}K_s \neq \emptyset.$ In fact, as $e^X \in \hbox{int}\mathcal{S}$ and $[X_s,X_n]=0$, we get $e^{X_s}e^{X_n} \in \hbox{int}\mathcal{S}$. By the continuity of product, there is open set $V$ containing $e^{X_s}$ such that $V e^{X_n} \subset \hbox{int}\mathcal{S}$. Then $V \in \hbox{int}K_s$. It is not hard to prove that $K_s$ is a semigroup of $K$, which implies that $K_s$ contains the identity component of $K$. Therefore $\exp(\tau X_n) \in \hbox{int}\mathcal{S}$, for some $\tau \in \R$ with $\hbox{ad}{(\tau X_n)}$ nilpotent.  By the Lemma (\ref{lemma}),  $S = \mathcal{S}$. 

Consider the set 
\begin{equation*}
    \mathfrak{h} = \{W \in \mathfrak{g}^0: \mathcal{D}(W) \in \mathfrak{r}(\mathfrak{g}^0)\}. 
\end{equation*}

We claim that $\mathfrak{h}$ is a $df_0-$invariant subalgebra of $\mathfrak{g}^0$ and $Z \in \mathfrak{h}$. In fact, by construction $\mathfrak{h}$ is a subspace of $\mathfrak{g}^0$. The Lie bracket stability follows by the Jacobi identity and the fact of $\mathfrak{r}(\mathfrak{g}^0)$ is an ideal of $\mathfrak{g}^0$. Since $\mathfrak{r}(\mathfrak{g}^0)$ is $\mathcal{D}-$invariant, we get that $\mathfrak{h}$ is $\mathcal{D}-$invariant such that $\mathcal{D}(\mathfrak{h}) \subset \mathfrak{r}(\mathfrak{g}^0)$. 

Also, as $R$ is a $f_0$-invariant solvable Lie group of $G^0$, by Proposition (\ref{lemmaimport}) we obtain $R \subset \mathcal{R}$. Considering $H$ as the connected Lie subgroup of $G^0$ generated by $\mathfrak{h}$, by the proposition (\ref{lemmaimport}) we obtain $H \subset \mathcal{R}$. By the lemma (\ref{AginA}) we have 
\begin{equation*}
    \mathcal{R}_k\cdot \exp(kZ) \subset \mathcal{R}\cdot \exp(kZ) \subset \mathcal{R}. 
\end{equation*}

In particular, it follows that $\mathcal{S} \subset \pi(G^0 \cap \mathcal{R})$. Then 
\begin{equation*}
    G^0 / R \subset \pi(\mathcal{R} \cap G^0). 
\end{equation*}

Therefore 
\begin{equation*}
    G^0 \subset (\mathcal{R} \cap G^0)R \subset \mathcal{R} R \subset \mathcal{R}. 
\end{equation*}
\QEDA

We now are able to prove the main result associated with the system $\Sigma^C$. 

\begin{proposition} Let $G$ be a connected Lie group with finite semisimple center and consider the linear system (\ref{linsyst}) defined over $G$. If $e \in \hbox{int}\mathcal{R}$ then $G^{0,+} \subset \mathcal{R}$.   
\end{proposition}

\noindent\textit{Proof:} For any $g \in G^{+}$, there is a $k \in \N$ such that $f_0^{-k}(g) \in \mathcal{R}$, since $\mathcal{R}$ is open and $G^+$ is stable in negative time. Then $g \in f_0^k\left(\mathcal{R}\right) \subset \mathcal{R}$. The previous proposition ensure that $G^0 \subset \mathcal{R}$. Using the lemma (\ref{AginA}), we have $G^{0,+} \subset \mathcal{R}$. \QEDA

By considering the same reasoning for the reversed-time system, we get by the theorem above that $G^{0,-} \subset \mathcal{C}. $ Finally, we have the following condition for controllability.

\begin{theorem}\label{finitecenter}Let $G$ be a Lie groups with finite semisimple center. If $\mathcal{R}$ is open and every eigenvalue of $\mathcal{D}$ has zero real part, then the system $(\Sigma)$ is controllable.
\end{theorem}

\noindent\textit{Proof:} As a matter of fact, every eigenvalue of $\mathcal{D}$ has zero real part if, and only if, $df_0$ has only eigenvalues with norm $1$. Therefore $G = G^0$. As $\mathcal{R}$ is open, we obtain 
\begin{equation*}
    G^0 = \mathcal{R} \cap \mathcal{C}, 
\end{equation*}
which implies in controllability. \QEDA 

Now, for the function 
\begin{equation*}
    d\hat{f}_0 = e^{\mathcal{D}_1} \circ ... \circ e^{\mathcal{D}_n}, 
\end{equation*}
for $\mathcal{D}_j \in \hbox{Der}(\mathfrak{g})$ for every $j=1,...,n$ and $\mathcal{D}_j \circ \mathcal{D}_k = \mathcal{D}_k \circ \mathcal{D}_j$, for $j,k=1,...,n$. Then $e^{\mathcal{D}_j}e^{\mathcal{D}_k} =e^{\mathcal{D}_k}e^{\mathcal{D}_j}$ for every $k,j=1.,...,n$. By induction hypothesis, suppose that the following statement holds for the automorphism $d\hat{f}_0$.  

\begin{claim}\label{claim1}Let $G$ be a Lie groups with finite semisimple center. If $\mathcal{D}_j \circ \mathcal{D}_k = \mathcal{D}_k \circ \mathcal{D}_j$, for $j,k=1,...,n$, $\mathcal{R}$ is open and every eigenvalue of $\mathcal{D}_j$ has zero real part, then the system $(\Sigma)$ is controllable.
\end{claim}

Now, take $df_0 = e^{\mathcal{D}_1} \circ ... \circ e^{\mathcal{D}_n}  \circ e ^{\mathcal{D}_{n+1}}$, that is, 
\begin{equation*}
    df_0 = d\hat{f}_0 \circ e^{\mathcal{D}_{n+1}}. 
\end{equation*}

Let us suppose that all eigenvalues of $\mathcal{D}_j$ has zero real part for every $j=1,...,n,n+1$ and $\mathcal{R}$ is open. As all the exponentials above commutes, the eigenvalues of $df_0$ are the product of all eigenvalues of $e^{\mathcal{D}_j}, j=1,...,n$ (see the Remark (\ref{eigenvalues})). Also, considering $\mathfrak{g}^0$ and  $\hat{\mathfrak{g}}^0$ as the central Lie algebras of $df_0$ and $d\hat{f}_0$ respectivelly $G^0$ and $\hat{G}^0$ its central Lie subgroups. We get that 
$\mathfrak{g}^0 = \hat{\mathfrak{g}}^0$. Consequently, by the claim above, we get 
\begin{equation*}
    G = \hat{G}^0 = G^0 = \mathcal{R} \cap \mathcal{C}
\end{equation*}
that is, the system is controllable. 

Therefore, considering $f_0$ satisfying 
\begin{equation*}
    df_0(X) = e^{\mathcal{D}_1} \circ ... \circ e^{\mathcal{D}_n}(X), X \in \mathfrak{g},
\end{equation*}
the next theorem summarize all the results proved so far in this chapter. 

\begin{theorem}Let $G$ be a Lie groups with finite semisimple center. If $\mathcal{D}_j \circ \mathcal{D}_k = \mathcal{D}_k \circ \mathcal{D}_j$, for $j,k=1,...,n$, $\mathcal{R}$ is open and every eigenvalue of $\mathcal{D}_j$ has zero real part, then the system $(\Sigma)$ is controllable.
\end{theorem}

\begin{remark}\label{eigenvalues}Let $V$ be a finite dimensional vector space and $\mathcal{L}(V)$ be the set of all linear operators of $V$. If $A,B \in \mathcal{L}(V)$ are linear maps such that $AB = BA$, it follows by \cite[Exercise 20.10]{elonlin} that there is a basis of $V$ such that $A$ and $B$ are triangular superior. This means that there is a change of basis matrix $P$ such that $PAP^{-1}$ and $PBP^{-1}$ are also triangular superior. For the product $AB$, we get 
\begin{equation*}
    PABP^{-1} = PAP^{-1}PBP^{-1}. 
\end{equation*}
which is the product of two triangular superior matrices, that is also a triangular superior matrix. As the eigenvalues of $A$ and $B$ are the elements of the diagonal of $PAP^{-1}$ and $PBP^{-1}$ respectivelly, the eigenvalues of $PABP^{-1}$ are the product of the eigenvalues of $A$ and $B$.     
\end{remark}

\begin{remark}Following \cite{AyalaeAdriano2} let us define the environment of the continuous case. Consider the family of ODE's
\begin{equation}\label{ode}
    \dot{g}(t) = \mathcal{X}(g(t)) + \sum_{i=1}^m u_j(t)X^j(g(t)), 
\end{equation}
where $\mathcal{X}$ is a linear vector field on $G$, $X^j$ are right invariant vector fields on $G$ and $u \in \mathcal{U} \subset L^{\infty}(\R, \Omega \subset \R^m)$ is the class of admissible controls with $\Omega$ a convex subset of $\R^m$. Denoting by $\Phi_{t,u}(g)$ the flow of (\ref{ode}) and $\phi_t$ the flow of $\mathcal{X}$, it is known of the Lie theory that $\phi$ is a $1-$parameter group of automorphism of $G$. Taking the derivation $\mathcal{D}: \mathfrak{g} \longrightarrow \mathfrak{g}$ defined by $\mathcal{D}(Y) = -[\mathcal{X},Y](e)$, the relation between $\mathcal{D}$ and $\phi_t$ is given by
\begin{equation}
    d\phi_t = e^{t\mathcal{D}}, \forall t \in \R. 
\end{equation}

Considering $f_0(g) = \phi_1(g)$ and pointwise controls 
\begin{equation*}
    u=(...,u(-1),u(0),u(1),u(2),...) = (...,u_{-1},u_0,u_1,u_2,...) \in \mathcal{U}
\end{equation*}
with $u_k \in \Omega$, we can define the system
\begin{equation}\label{odediscrete}
    x_{k+1} = f_{u_k}(x_k), k \in \N_0,
\end{equation}
where $f: \Omega \times G \longrightarrow G$ is given by 
$f_{u_0}(g) = \Phi_{1,u}(g)$, it follows that the system (\ref{odediscrete}) is a discrete-time linear system on $G$. As a matter of fact, the solution of (\ref{ode}) satisfies 
\begin{equation*}
    \Phi_{\tau, u}(g) = \Phi_{\tau, u}(e) \phi_t(g), \forall \tau \in \R. 
\end{equation*}

Then $f_{u_0}(g) = \Phi_{1,u}(g) = \Phi_{1,u}(e) \phi_t(g) = f_{u_0}(e)f_0(g)$ for every $u \in \mathcal{U}$, with $f_0: G \longrightarrow G$ an automorphism by construction. According to the expression (5) on \cite{AyalaeAdriano2}, there is a $W \in \mathfrak{g}^0$ such that the function $f_0$ can be considered as $f_0^k(g) = e^{(kW)} g e^{(-kW)} h$, with $h = h_{W,k,g} \in R = \langle \exp(\mathfrak{r}(\mathfrak{g}^0)) \rangle$ and $k \in \N$, which is preciselly the case in the theorem (\ref{finitecenter}). In this particular case, very hypothesis we supposed to be valid here is fulfilled in the continuous case.
\end{remark}

\begin{example} Let us consider the Lie group $G = SL_2(\R)$, given by 
\begin{equation*}
    \hbox{SL}_2(\R): \left\{\begin{bmatrix}
    a & b\\
    c & d
    \end{bmatrix} \in \hbox{GL}_2(\R): ad - cb = 1\right\}. 
\end{equation*}

This group is a connected Lie subgroup of $GL_2(\R)$ with Lie algebra 
\begin{equation*}
    \mathfrak{sl}_2(\R) = \left\{\begin{bmatrix}
    a & b\\
    c & d
    \end{bmatrix} \in \mathfrak{gl}_2(\R): a + d = 0\right\}.
\end{equation*}

In particular, the Lie algebra $\mathfrak{sl}_2(\R) $ is semisimple and, consequently, the Lie group $\hbox{SL}_2(\R)$ is a semisimple connected Lie group \cite[Definition 1.30]{sanmartin2}. It is known that $\hbox{Aut}(\mathfrak{sl}_2(\R)) = \hbox{Inn}(\mathfrak{sl}_2(\R))$, that is, every automorphism of $\mathfrak{sl}_2(\R)$ is inner in the sense of if $T \in \hbox{Aut}(\mathfrak{sl}_2(\R))$, there are $Y_1,...,Y_n \in \mathfrak{gl}_2(\R)$ such that 
\begin{equation*}
    T(X) = e^{\adj{Y_1}}...e^{\adj{Y_{n}}}(X). 
\end{equation*}

Using \cite[Proposition 5.15]{sanmartin1}, this immediately implies that, considering $h = e^{Y_1}...e^{Y_n} \in GL_{2}(\R)$, we have that the conjugation $C_h(g) = h g h^{-1}$ has as differential at the identity the function $T$ above. This allows us to define the class of linear systems of $\hbox{SL}_2(\R)$. In fact, given a $h \in \hbox{GL}_2(\R)$, consider a function $f: U \times \hbox{SL}_2(\R) \longrightarrow \hbox{SL}_2(\R)$ given by 
\begin{equation*}
    f_u(g) = 
    \begin{bmatrix}
        f^{11}_u(e) & f^{12}_u(e)\\
        f^{21}_u(e) & f^{22}_u(e)
    \end{bmatrix} hgh^{-1}
\end{equation*}
such that $f^{11}_0(e) = f^{22}_0(e) = 1$, $f^{21}_0(e) = f^{12}_0(e) = 0$ and $f^{11}_u(e)f^{22}_u(e) - f^{21}_u(e)f^{12}_u(e)=1$ for all $u \in U$. Considering the discrete-time system
\begin{equation*}
    \Sigma_L: x_{k+1} = f_{u_k}(x_k), k \in \N_0,
\end{equation*}
we claim that the system $(\Sigma_L)$ is a linear system on $\hbox{SL}_2(\R)$. In fact, given the properties of the matrix $f_u(e)$, we have $f_0(g) = hgh^{-1}$. Then $f_0$ is an automorphism of $\hbox{SL}_2(\R)$. Besides, by construction we have $f_u(g) = f_u(e)f_0(g)$. In particular, every linear system in $SL_2(\R)$ has the form of $(\Sigma_L)$ given the product property in the definition (\ref{linearsystem}) of linear systems. 

Now, take the matrices $h$ and $h^{-1}$ as
\begin{equation*}
    h= 
    \begin{bmatrix}
    h_{11} & h_{12}\\
    h_{21} & h_{22}
    \end{bmatrix} \hbox{ and }
    h^{-1} = \frac{1}{h_{11}h_{22} - h_{21}h_{12}}
    \begin{bmatrix}
    h_{22} & -h_{21}\\
    -h_{12} & h_{11}
    \end{bmatrix}, 
\end{equation*}
and a element $g \in \hbox{SL}_2(\R)$ in the form $g = \begin{bmatrix}
    g_{11} & g_{12}\\
    g_{21} & g_{22}
    \end{bmatrix}$. We have 
\begin{equation*}
    f_u(g) = 
    \begin{bmatrix}
        f^{11}_u(e) & f^{12}_u(e)\\
        f^{21}_u(e) & f^{22}_u(e)
    \end{bmatrix} 
    f_0(g),
\end{equation*}
where $f_0(g) = h g h^{-1}$ is given by
\begin{equation*}
     f_0(g) = 
    \begin{bmatrix}
       \frac{g_{12} h_{11} h_{21} + g_{22} h_{12} h_{21} - g_{11} h_{11} h_{22} - g_{21} h_{12} h_{22}}{h_{12} h_{21} - h_{11} h_{22}} & \frac{(g_{12} h_{11}^2 - h_{12} (g_{11} h_{11} - g_{22} h_{11} + g_{21} h_{12}))}{(-h_{12} h_{21} + h_{11} h_{22})}\\ 
       \frac{(-g_{12} h_{21}^2 + h_{22} (g_{11} h_{21} - g_{22} h_{21} + g_{21} h_{22}))}{(-h_{12} h_{21} + h_{11} h_{22})} & \frac{(g_{12} h_{11} h_{21} + g_{22} h_{11} h_{22} - h_{12} (g_{11} h_{21} + g_{21} h_{22}))}{(-h_{12} h_{21} +h_{11} h_{22})}
    \end{bmatrix}
\end{equation*}  

In particular, the center of $\hbox{SL} _2(\R)$ is center is $\Z_2$. Then we can apply our results to this set. We claim that, for every $h \in \hbox{GL}_2(\R)$, the linear map $df_0$ always have an eigenvalue $\lambda = 1$. In fact, considering the function $f_0$ as a function of $f_0: \R^4 \longrightarrow \R^4$, we have that  
\begin{equation*}
    df_0 = 
    \left[
\begin{array}{cccc}
    -\frac{h_{11} h_{22}}{h_{12} h_{21}-h_{11} h_{22}} & \frac{h_{11} h_{21}}{h_{12} h_{21}-h_{11} h_{22}} & -\frac{h_{12} h_{22}}{h_{12} h_{21}-h_{11} h_{22}} & \frac{h_{12} h_{21}}{h_{12} h_{21}-h_{11} h_{22}} \\
    -\frac{h_{11} h_{12}}{h_{11} h_{22}-h_{12} h_{21}} & \frac{h_{11}^2}{h_{11} h_{22}-h_{12} h_{21}} & -\frac{h_{12}^2}{h_{11} h_{22}-h_{12} h_{21}} & \frac{h_{11} h_{12}}{h_{11} h_{22}-h_{12} h_{21}} \\
    \frac{h_{21} h_{22}}{h_{11} h_{22}-h_{12} h_{21}} & -\frac{h_{21}^2}{h_{11} h_{22}-h_{12} h_{21}} & \frac{h_{22}^2}{h_{11} h_{22}-h_{12} h_{21}} & -\frac{h_{21} h_{22}}{h_{11} h_{22}-h_{12} h_{21}} \\
    -\frac{h_{12} h_{21}}{h_{11} h_{22}-h_{12} h_{21}} & \frac{h_{11} h_{21}}{h_{11} h_{22}-h_{12} h_{21}} & -\frac{h_{12} h_{22}}{h_{11} h_{22}-h_{12} h_{21}} & \frac{h_{11} h_{22}}{h_{11} h_{22}-h_{12} h_{21}} \\
\end{array}\right]
\end{equation*}

The characteristic polynomial of $df_0$ is given by 
\begin{equation*}
    p(\lambda) = \frac{(\lambda -1)^2 \left((\lambda  h_{22}-h_{11}) (h_{22}-\lambda  h_{11})+(\lambda +1)^2 h_{12} h_{21}\right)}{h_{12} h_{21}-h_{11} h_{22}}.
\end{equation*}

Hence, the spectrum of $df_0$ is the set
\begin{equation*}
    \hbox{Spec}(df_0) = \left\{1,\lambda_1,\lambda_2\right\}
\end{equation*}
where 
\begin{eqnarray*}
    \lambda_1 &=& \frac{-(h_{11}+h_{22}) \sqrt{h_{11}^2-2 h_{11} h_{22}+4 h_{12} h_{21}+h_{22}^2}+h_{11}^2+2 h_{12} h_{21}+h_{22}^2}{2 (h_{11} h_{22}-h_{12} h_{21})}, \\
    \lambda_2 &=&\frac{(h_{11}+h_{22}) \sqrt{h_{11}^2-2 h_{11} h_{22}+4 h_{12} h_{21}+h_{22}^2}+h_{11}^2+2 h_{12} h_{21}+h_{22}^2}{2 (h_{11} h_{22}-h_{12} h_{21})}.
\end{eqnarray*}

In particular, $df_0$ has real spectrum if, and only if $h_{11}^2-2 h_{11} h_{22}+4 h_{12} h_{21}+h_{22}^2 \geq 0$. The algebraic multiplicity of $1$ is two. 

Let us consider the case when $U \subset \R$. As the dimension of the group is $3$, the minimum possible time that makes the regular reachable set $\hat{\mathcal{R}}_k$ to contain $e$ is $k = 3$, given that the matrix $\frac{\partial }{\partial (u,v)} f_u \circ f_v(e)$ has the form
\begin{equation*}
    \frac{\partial }{\partial (u,v)} f_u \circ f_v(e) = 
    \left[ 
    \begin{array}{cc}
         v_1& v_2  \\
         v_3&v_4 
    \end{array}
    \right]
\end{equation*}
with $v_1,v_2,v_3,v_4 \in \R^2$. For higher dimensions, the reasoning would be the same. 

Let us explore some numeric examples. Take the case when $U \subset \R$ is a compact convex neighborhood of $0$ and 
\begin{equation}\label{matrixh}
    h = 
    \begin{bmatrix}
        1&1\\
        0&1
    \end{bmatrix}, 
\end{equation}
with 
\begin{equation*}
    f_u(e) = 
    \begin{bmatrix}
        1+u&-u\\
        u&1-u
    \end{bmatrix}.
\end{equation*}

Considering the $2-$step nilpotent matrix
\begin{equation*}
    M = \begin{bmatrix}
        0 & 1 \\
        0 & 0
    \end{bmatrix}. 
\end{equation*}
we get that 
\begin{equation}
    e^M = \sum_{ n \in \N} \frac{M^n}{n!} = I + M = 
    \begin{bmatrix}
        1&1\\
        0&1
    \end{bmatrix} = h 
\end{equation}
At first, let us check the eigenvalues of the function $df_0$ with $h$ defined above. In fact, considering  
\begin{equation*}
    f_0(g_{11},g_{12}, g_{21}, g_{22}) = \left(g_{11}+g_{21} ,-g_{11}+g_{12}-g_{21}+g_{22}, g_{21}, g_{22}-g_{21}
    \right)
\end{equation*}
with derivative given by 
\begin{equation*}
    df_0 = \left[
    \begin{array}{cccc}
        1 & 0 & 1 & 0 \\
        -1 & 1 & -1 & 1 \\
        0 & 0 & 1 & 0 \\
        0 & 0 & -1 & 1 \\
    \end{array}\right]
\end{equation*} 
the spectrum of $df_0$ is $\hbox{Spec}(df_0) = \{1\}$. Let us check the openness of $\mathcal{R}$. It is simple to verify that $f_u(e) \in \hbox{SL}_2(\R),$ for every $u \in U$. Following the definition (\ref{regular}), we claim that $e \in \hat{\mathcal{R}}$. In fact, take $(u,v,w) \in \hbox{int}U^3$. For $k=3$, using the notation $f_{u,v,w} = f_{u} \circ f_{v} \circ f_{w}$, we have 
\begin{equation*}
    f_{u,v,w}(e) = f_{u}(e) C_h (f_{v}(e)) C_{h^2}(f_{w}(e))
\end{equation*}

The matrix above is given by

\resizebox{\linewidth}{!}{%
$\displaystyle
    f_{u,v,w}(e)=\left[    
    \begin{array}{cc} 
        (w+1)((u+1) (2 v+1)-u v)&  (-u (1 - 2 v) - 4 (1 + u) v) + (-u v + (1 + u) (1 + 2 v))\\
        ((1 - u) v + u (1 + 2 v)) (1 + w) & ((1 - u) (1 - 2 v) - 4 u v) + ((1 - u) v + u (1 + 2 v))
    \end{array} \right]
    $}
where $f_{u,v,w} = f_{u} \circ f_{v} \circ f_{w}$, whose derivative is given by 
\begin{equation}
    \frac{\partial }{\partial(u,v,w)} f_{u} \circ f_{v} \circ f_{w}(e) =  
    \begin{bmatrix}
        \begin{bmatrix}
            (1+v)(1+w)\\
            (2 + u) (1 + w)\\
            1+ 2v + u (1+ v)
        \end{bmatrix} &
        \begin{bmatrix}
            -v\\ -2-u\\
            0
        \end{bmatrix}\\
        \begin{bmatrix}
            (1+v)(1+w)\\ 
            (1+u)(1+w)\\ 
            v + u (1+v)
        \end{bmatrix} &
        \begin{bmatrix}
           -v\\ -1-u\\ 0
        \end{bmatrix}
    \end{bmatrix}.
\end{equation}

Taking the vectors of the matrix above, one can prove that the subespace generated by them is $3-$dimensional. Hence, the matrix above has rank $3$, for every $u \in \hbox{int}U$. Then $e \in \hat{\mathcal{R}}_3 \subset \hat{\mathcal{R}}$. As the set $\hat{\mathcal{R}}$ is open and $\hat{\mathcal{R}} \subset \mathcal{R}$, we have $e \in \hbox{int}\mathcal{R}$. By the Proposition (\ref{reachablesetprop}), item 6, the set $\mathcal{R}$ is open. This implies by the Proposition (\ref{openess}) that $\hbox{int}\mathcal{C}$ is open. Therefore, by the Theorem (\ref{finitecenter}), the system is controllable.  
\end{example}


\begin{thebibliography}{99}

\bibitem{AyalaeAdriano2}AYALA, V.; SILVA, A. \emph{ Controllability of linear systems on Lie groups with finite
semisimple center.} SIAM J. control and optimization, v. 55, (2017).

\bibitem{AyalaeRomaneAdriano}AYALA, V.; SILVA, A.; ROMAN-FLORES, H. \emph{ The dynamic of a Lie group
endomorphism.} Open Math, v. 15, (2017).

\bibitem{TAJ} CAVALHEIRO, T.; SANTANA, A. J.; COSSICH, J. A. N. \emph{ Controllability of discrete-time linear systems on solvable Lie groups.} arXiv preprint arXiv:2302.00145, (2023).

\bibitem{CCS1}COLONIUS, F.; SANTANA, A. J.; COSSICH, J. \emph{ Outer invariance entropy for discrete-time linear systems on Lie groups.} ESAIM: Control, Optimisation and Calculus of Variations, v. 27, (2021).


\bibitem{CCS2}COLONIUS, F.; SANTANA, A. J.; COSSICH, J. \emph{ Controllability properties and invariance pressure for discrete-time linear systems .} Journal of Dynamics and Differential Equations, v. 34, (2022).


\bibitem{elonlin} LIMA, E. L. \emph{ Algebra Linear.} IMPA. (2020) 

\bibitem{fresnel}FRESNEL, J. \emph{ Compact subgroups of $\hbox{GL}(n,\mathbb{C})$.} Rend. Sem. Mat. Univ. Pandova. (2006)

\bibitem{jouan}JOUAN, P. \emph{ Controllability of Linear Systems on Lie Groups.} J. Dyn. Control Syst.,
v. 17, (2011).

\bibitem{helgason} HELGASON. S. \emph{ Differential geometry, Lie groups and Symmetric spaces.} Academic Press. (1978)

\bibitem{pontr} PONTRJAGIN, L. \emph{ Topological Groups.} PRinceton University Press. (1946)

\bibitem{onish}ONISHCHIK, A.; VINBERG, E. Lie groups and Lie algebras. Springer, (1993). v. 41.

\bibitem{sanmartin2}SAN-MARTIN, L. A. B. \emph{ Algebras de Lie.} Editora Unicamp, 2010.

\bibitem{sanmartin1}SAN-MARTIN, L. A. B. \emph{ Lie groups.} Springer, (2016).

\bibitem{sanmartin3} SAN-MARTIN, L. A. B. \emph{ Invariant control sets on flag manifolds.} Math. Control Signals Systems. vol. 6. (1993) pg 41-61. 

\bibitem{sontag1}SONTAG, E. \emph{ Mathematical Control Theory. Deterministic finite- dimensional
systems.} Springer-Verlag, (1998).

\bibitem{wustner} WUSTNER, M. \emph{ On the surjectivity of the exponential function on solvable Lie groups.} Math. Nachr. (1998), pg 255 - 266.

\end{thebibliography}
\end{document}